\magnification=\magstep1
\input amstex
\input amsppt.sty
\hsize=5.5truein
\vsize=8.5truein

\define\a{\alpha}
\define\be{\beta}
\define\ga{\gamma}
\define\de{\delta}

\define\n{\Bbb N}
\define\lo{\Bbb L}
\define\p{\Cal P}
\define\q{\Bbb Q}
\define\z{\Bbb Z}

\define\mm{\Cal M}
\define\uu{\Cal U}
\define\g{\Cal G}
\define\rvol{\operatornamewithlimits{rvol}}

\define\Vol{\operatornamewithlimits{Vol}}
\define\ar{\operatornamewithlimits{Area}}

\NoRunningHeads
\NoBlackBoxes
\parskip=9pt
\topmatter
\title{\bf UNIVERSAL COUNTING OF LATTICE POINTS IN POLYTOPES}
\endtitle
\endtopmatter

\flushpar
\hphantom{xxxxxxxxx}{\smc Imre B\'ar\'any} \footnote{ Partially
supported by Hungarian Science Foundation
Grant  T 016391, and by the French--Hungarian Bilateral
Project Balaton F--30/96}
\newline
\hphantom{xxxxxxxxx}Mathematical Institute
of the Hungarian Academy of Sciences,
\newline
\hphantom{xxxxxxxxx}POB 127, 1364 Budapest, Hungary
\newline
\hphantom{xxxxxxxxx}{\tt barany\@math-inst.hu}, and
\newline
\hphantom{xxxxxxxxx}Dept. Mathematics,
\newline
\hphantom{xxxxxxxxx}University College London,
\newline
\hphantom{xxxxxxxxx}Gower Street, WC1E 6BT London, UK
\flushpar
\hphantom{xxxxxxxxx}{\smc Jean--Michel Kantor}
\newline
\hphantom{xxxxxxxxx}Institut Math\'ematique de Jussieu,
\newline
\hphantom{xxxxxxxxx}Universit\'e Paris 7, 4 Place Jussieu,
\newline
\hphantom{xxxxxxxxx}75252 Paris, France
\newline
\hphantom{xxxxxxxxx}{\tt kantor\@math.jussieu.fr}

\document

\bigskip
{\bf Abstract.} Given a lattice polytope $P$ (with
underlying lattice $\lo$), the universal counting function
$\uu_P(\lo')=|P\cap \lo'|$ is defined on all lattices $\lo'$
containing $\lo$. Motivated by questions concerning lattice
polytopes and the Ehrhart polynomial, we study the equation
$\uu_P=\uu_Q$.

\baselineskip=15.4pt
\vskip1cm
Mathematics Subject Classification: 52B20, 52A27, 11P21

\pagebreak

\bigskip
{\smc 1. The universal counting function}

We will denote by $V$ a vector space of dimension
$n$, by $\lo$ a lattice in $V$, of rank $n$. Let
$$
\g_{\lo}=\lo \rtimes GL(\lo)
$$
be the group of affine maps of $V$ inducing isomorphism of
$V$ and $\lo$ into itself; in case
$$
\lo=\z^n \subset V=\q^n, \g_n =\z^n \rtimes GL(\z^n)
$$
corresponds to affine unimodular maps. An $\lo$--polytope
is the convex hull of finitely many points from $\lo$;
$\p_{\lo}$ denotes the set of all $\lo$--polytopes. For
a finite set $A$ denote by $|A|$ its cardinality.
Finally, let $\mm_{\lo}$ be the set of all lattices
containing $\lo$.

\demo{Definition 1} Given any $\lo$--polytope $P$, the
function $\uu_P:\mm_{\lo} \to \z$ defined by
$$
\uu_P(\lo')=|P\cap \lo'|
$$
is called the {\it universal counting function} of $P$.
\enddemo

This is just the restriction of another function $\uu:
\p_{\lo}\times \mm_{\lo} \to \z$ to a fixed $P\in \p_{\lo}$, where
$\uu$ is given by
$$
\uu(P,\lo')=|P\cap \lo'|.
$$
Note, further, that $\uu_P$ is invariant under the group,
$\g_{tr}$, generated by $\lo$--translations and the reflection
with respect to the origin, but, of course, not invariant
under $\g_{\lo}$.

Example 1. Take for $\lo'$ the lattices
$\lo_k=\frac 1k \lo$ with $k\in \n$. Then
$$
\uu_P(\lo_k)=|P\cap \frac 1k \lo|=|kP\cap \lo|=E_P(k)
$$
where $E_P$ is the Ehrhart polynomial of $P$ (see [Ehr]). We
will need some of its properties that are described in the
following theorem (see for instance [Ehr],[GW]).
Just one more piece of notation:
if $F$ is a facet of $P$ and $H$ is the affine hull of $F$,
then the relative volume volume of $F$ is defined as
$$
\rvol (F)=\frac {\Vol_{n-1} (F)}{\Vol_{n-1} (D)}
$$
where $D$ is the fundamental parallelotope of the
$(n-1)$--dimensional sublattice of $H \cap \lo$. For a face $F$ of $P$
that is at most $(n-2)$--dimensional let $\rvol (F)=0$.
Note that the relative volume is
invariant under $\g_{\lo}$ and can be computed,
when $\lo=\z^n$, since then the denominator is the euclidean
length of the (unique) primitive outer normal to $F$ (when $F$ is a facet).

\proclaim{Theorem 1} Assume $P$ is an $n$--dimensional
$\lo$--polytope. Then $E_P$ is a polynomial in $k$ of
degree $n$. Its main coefficient is $\Vol (P)$,
and its second coefficient equals
$$
\frac 12 \sum _{F \text{ a facet of } P} \rvol (F).
$$
\endproclaim

It is also known that $E_P$ is a $\g_{\lo}$--invariant
valuation, (for the definitions see [GW] or [McM]).
The importance of $E_P$ is reflected in the
following statement from [BK]. For a
$\g_{\lo}$--invariant valuation $\phi$
from $\p_{\lo}$ to an abelian group $G$, there exists a
unique $\ga=(\ga_i)_{i=0,\dots,n}$ with $\ga_i \in G$ such that
$$
\phi(P)=\sum \ga_ie_{P,i}
$$
where $e_{P,i}$ is the coefficient of $k^i$ of the Ehrhart
polynomial.

It is known that $E_P$ does not determine $P$,
even within $\g_{\lo}$ equivalence. [Ka] gives examples
lattice--free $\lo$--simplices with identical Ehrhart
polynomial that are different under $\g_{\lo}$. The aim of this
paper is to investigate whether and to what extent the
universal counting function determines $P$.

We give another description of $\uu_P$. Let $\pi\:V \to V$ be
any isomorphism satisfying $\pi(\lo) \subset \lo$. Define, with
a slight abuse of notation,
$$
\uu_P(\pi)=|\pi(P)\cap \lo|=|P\cap \pi^{-1}(\lo)|.
$$
Set $\lo'=\pi^{-1}(\lo)$. Since $\lo'$ is a lattice containing
$\lo$ we clearly have
$$
\uu_P(\pi)=\uu_P(\lo').
$$
Conversely, given a lattice $\lo' \in \mm_{\lo}$, there is an
isomorphism $\pi$ satisfying the last equality. (Any linear
$\pi$ mapping a basis of $\lo$ to a basis of $\lo'$ suffices.)
The two definitions of $\uu_P$ via lattices or isomorphisms with
$\pi(\lo) \subset \lo$ are equivalent. We will use the common
notation $\uu_P$.

Example 2. Anisotropic dilatations. Take $\pi :\z^n \to
\z^n$ defined by
$$
\pi(x_1,\dots,x_n)=(k_1x_1,\dots,k_nx_n),
$$
where $k_1,\dots,k_n \in \n$. The corresponding map $\uu_P$
extends the notion of Ehrhart polynomial and Example 1.

Simple examples show that $\uu_P$ is not a polynomial in
the variables $k_i$.

\bigskip
{\smc 2. A necessary condition}

Given a nonzero $z \in \lo^*$, the dual of $\lo$, and an
$\lo$--polytope $P$, define $P(z)$ as the set of points in
$P$ where the functional $z$ takes its maximal value. As is
well known, $P(z)$ is a face of $P$. Denote by $H(z)$ the
hyperplane $z\cdot x=0$ (scalar product). $H(z)$ is clearly
a lattice subspace.

\proclaim{Theorem 2} Assume $P,Q$ are $\lo$--polytopes with
identical universal counting function. Then, for every primitive
$z \in \lo^*$,
$$
\rvol P(z) +\rvol P(-z)=\rvol Q(z) + \rvol Q(-z).
\tag *
$$
\endproclaim

The theorem shows, in particular, that if $P(z)$ or $P(-z)$
is a facet of $P$, then $Q(z)$ or $Q(-z)$ is a facet of
$Q$. Further, given an $\lo$--polytope $P$, there
are only finitely many possibilities for the outer normals
and volumes of the facets of another polytope $Q$ with
$\uu_P=\uu_Q$. So a well--known theorem of Minkowski implies,

\proclaim{Corollary 1} Assume $P$ is an $\lo$--polytope.
Then, apart from lattice translates, there are only
finitely many $\lo$--polytopes with the same universal
counting functions as $P$.
\endproclaim

\demo{Proof of Theorem 2} We assume that $P,Q$ are full--dimensional
polytopes. It is enough to prove the theorem in the special
case when $\lo =\z^n$ and $z=(1,0,\dots,0)$. There is nothing to
prove when none of $P(z),P(-z),Q(z)$, $Q(-z)$ is a facet since
then both sides of (*) are equal to zero. So
assume that, say, $P(z)$ is a facet, that is, $\rvol
P(z)>0$.

For a positive integer $k$ define the linear map $\pi_k \:V \to V$ by
$$
\pi_k(x_1,\dots,x_n)=(x_1,kx_2,\dots,kx_n).
$$
The
condition implies that the lattice polytopes $\pi_k(P)$ and
$\pi_k(Q)$ have the same Ehrhart polynomial. Comparing
their second coefficients we get,
$$
\sum _{F \text{ a facet of }P} \rvol \pi_k(F)
=\sum _{G \text{ a facet of }Q} \rvol \pi_k(G),
$$
since the facets of $\pi_k(P)$ are of the form
$\pi_k(F)$ where $F$ is a facet of $P$.

Let $\zeta=(\zeta_1,\dots,\zeta_n) \in
\z^{n*}$ be the (unique) primitive outer normal to the facet
$F$ of $P$. Then $\zeta'=(k\zeta_1,\zeta_2,\dots,\zeta_n)$
is an outer normal to $\pi_k(F)$, and so it is a positive
integral multiple of the unique primitive outer normal
$\zeta ''$, that is $\zeta'=m\zeta''$ with $m$ a positive
integer. When $k$ is a large
prime and $\zeta$ is different from $z$ and $\zeta_1\ne 0$,
then $m=1$ and $\rvol \pi_k(F) =O(k^{n-2})$. When
$\zeta_1=0$, then $m=1$, again, and the ordinary $(n-1)$--volume
of $\pi_k(F)$ is $O(k^{n-2})$. Finally, when $\zeta=\pm z$,
$\Vol \pi_k(F)=k^{n-1}\Vol F$.

So the dominant term, when $k \to \infty$, is
$k^{n-1}(\rvol P(z)+\rvol P(-z))$ since by our assumption
$\rvol P(z) >0$. \qed
\enddemo

\bigskip
{\smc 3. Dimension two}

Let $P$ be an $\lo$--polygon in $V$ of dimension two.
Simple examples show again that $\uu_P$ is not a polynomial
in the coefficients of $\pi$.

In the planar case we abbreviate $\rvol P(z)$ as $|P(z)|$.
Extending (and specializing) Theorem 1 we prove

\proclaim{Proposition 3} Suppose $P$ and $Q$ are
$\lo$--polygons. Then $\uu_P=\uu_Q$ if and only if the
following two conditions are satisfied:
\item{(i)} $\ar(P)=\ar(Q)$,
\item{(ii)} $|P(z)|+|P(-z)|=|Q(z)|+|Q(-z)|$ for every primitive
$z \in \lo^*$.
\endproclaim

\demo{Proof} The conditions are sufficient: (i) and (ii)
imply that, for any $\pi$, $\ar(\pi(P))=\ar(\pi(Q))$ and
$|\pi(P)(z)|+|\pi(P)(-z)|=|\pi(Q)(z)|+|\pi(Q)(-z)|$.
We use Pick's formula for $\pi(P)$, (see [GW], say):
$$
|\pi(P)\cup \lo|=\ar \pi(P)+\frac 12 \sum_{z\text{ primitive}}|\pi(P)(z)|+1.
$$
This shows that $\uu_P=\uu_Q$, indeed.

The necessity of (i) follows from Theorem 1 immediatley, (via the main
coefficient of $E_P$), and the necessity of (ii) is the
content of Theorem 2.\qed
\enddemo

\proclaim{Corollary 2} Under the conditions of Proposition 3
the lattice width of $P$ and $Q$, in any direction $z \in
\lo^*$ are equal.
\endproclaim

\demo{Proof} The lattice width, $w(z,P)$, of $P$ in
direction $z \in \lo^*$ is, by definition (see [KL],[Lo]),
$$
w(z,P)=\max \{z\cdot (x-y)\: x,y \in P\}.
$$
In the plane one can compute the width along the boundary
of $P$ as well which gives
$$
w(z,P)=\frac 12 \sum_{e} |z\cdot e|
$$
where the sum is taken over all edges $e$ of $P$. This
proves the corollary.\qed
\enddemo

\proclaim{Theorem 3} Suppose $P$ and $Q$ are
$\lo$--polygons. Then $\uu_P=\uu_Q$ if and only if the
following two conditions are satisfied:
\item{(i)} $\ar(P)=\ar(Q)$,
\item{(ii)} there exist $\lo$--polygons $X$ and $Y$
such that $P$ resp. $Q$ is a lattice translate of $X+Y$ and
$X-Y$ (Minkowski addition).
\endproclaim

Remark. Here $X$ or $Y$ is allowed to be a segment or even
a single point. In the proof we will ignore translates and
simply write $P=X+Y$ and $Q=X-Y$.

\demo{Proof} Note that (ii) implies the second condition in
Proposition 3. So we only have to show the necessity of (ii).

Assume the contrary and let $P,Q$ be a counterexample to
the statement with the smallest possible number of edges.
We show first that for every (primitive) $z \in \lo^*$ at least
one of the sets $P(z),P(-z),Q(z),Q(-z)$ is a point.

If this were not the case, all four segments would contain
a translated copy of the shortest among them, which, when
translated to the origin, is of the form $[0,t]$. But then
$P=P'+[0,t]$ and $Q=Q'+[0,t]$ with $\lo$--polygons $P',Q'$.

We claim that $P',Q'$ satisfy conditions (i) and (ii) of Proposition 3. This
is obvious for (ii). For the areas we have that $\ar P-\ar P'$ equals
the area of the parallelogram with base $[0,t]$ and height $w(z,P)$.
The same applies to $\ar Q-\ar Q'$, but there the height is
$w(z,Q)$. Then Corollary 2 implies the claim.

So the universal counting functions of $P',Q'$  are identical.
But the number of edges
of $P'$ and $Q'$ is smaller than that of $P$ and $Q$.
Consequently there are polygons $X'$, $Y$ with $P'=X'+Y$,
and $Q'=X'-Y$. But then, with $X=X'+[0,t]$, $P=X+Y$ and
$Q=X-Y$, a contradiction.

Next, we define the polygons $X,Y$ by specifying their
edges. It is enough to specify the edges of $X$ and $Y$
that make up the edges $P(z),P(-z), Q(z)$, $Q(-z)$ in $X+Y$ and
$X-Y$. For this end we orient the edges of $P$ and $Q$
clockwise and set
$$
P(z)=[a_1,a_2],P(-z)=[b_1,b_2],Q(z)=[c_1,c_2],Q(-z)=[d_1,d_2]
$$
each of them in clockwise order. Then
$$
a_2-a_1=\a t, b_2-b_1=\be t, c_2-c_1=\ga t, d_2-d_1=\de t
$$
where $t$ is orthogonal to $z$ and $\a,\ga \geq 0$, $\be,
\de \leq 0$ and one of them equals $0$. Moreover, by
condition (ii) of Proposition 3, $\a-\be=\ga-\de$.

Here is the definition of the corresponding edges, $x,y$
of $X,Y$:
$$
\align
x=\a t,  y=\be t \text{ if } &\de=0,\\
x=\be t, y= \a t \text{ if } &\ga=0,\\
x=\ga t, y=-\de t \text{ if } &\be=0,\\
x=\de t, y=-\ga t \text{ if } &\a=0.
\endalign
$$
With this definition, $X+Y$ and $X-Y$ will have exactly the
edges needed.
We have to check yet that the sum of the $X$
edges (and the $Y$ edges) is zero, otherwise they won't
make up a polygon. But $\sum (x+y) =0$ since this is the
sum of the edges of $P$, and $\sum (x-y)=0$ since this is
the sum of the edges of $Q$. Summing these two equations
gives $\sum x=0$, subtracting them yields $\sum y=0$.
\qed\enddemo

\bigskip
{\smc 4. An example and a question}

Let $X$, resp. $Y$ be the triangle with vertices
$(0,0),(2,0),(1,1)$, and $(0,0), (1,1), (0,3)$. As it turns
out the areas of $P=X+Y$ and $Q=X-Y$ are equal. So Theorem
3 applies: $\uu_P=\uu_Q$. At the same time, $P$ and $Q$ are
not congruent as $P$ has six vertices while $Q$ has only five.

However, it is still possible that polygons with the same
universal counting function are equidecomposable. Precisely,
$P_1,\dots ,P_m$ is said to be a subdivision of $P$ if the $P_i$ are
$\lo$--polygons with pairwise relative interior, their
union is $P$, and the intersection of the closure of any
two of them is a face of both. Recall from section 1
the group $\g_{tr}$ generated by $\lo$--translations and
the reflection with respect to the origin. Two $\lo$--polygons $P,Q$ are
called $\g_{tr}$--equidecomposable if there are subdivisions
$P=P_1\cup \dots \cup P_m$ and $Q=Q_1\cup \dots \cup Q_m$
such that each $P_i$ is a translate, or the reflection of a
translate of $Q_i$ with the extra condition that $P_i$ is
contained in the boundary of $P$ if and only if $Q_i$ is
contained in the boundary of $Q$.

We finish the paper with a question which has connections
to a theorem of the late Peter Greenberg [Gr].
Assume $P$ and $Q$ have the same universal
counting function. Is it true then that they are
$\g_{tr}$--equidecomposable? In the example above, as in many other
examples, they are.

\bigpagebreak {\smc References}

\item{[BK]} {\smc U. Betke, M. Kneser}, Zerlegungen und
Bewertungen von Gitterpolytopen, {\it J. Reine ang. Math.} {\bf
358} (1985), 202--208.

\item{[Eh]} {\smc E. Ehrhart}, {\it Polinomes arithm\'etiques et
m\'etode des poly\'edres en combinatoire}, Birkhauser, 1977.

\item{[Gr]} {\smc P. Greenberg}, Piecewise
$SL_2$--geometry, {\it Transactions of the AMS}, {\bf 335}
(1993), 705--720.

\item{[GW]} {\smc P. Gritzmann, J. Wills}, Lattice points,
in: {\it Handbook of convex geometry, ed. P. M. Gruber, J.
Wills}, North Holland, Amsterdam, 1988.

\item{[KL]} {\smc R. Kannan, L. Lov\'asz}, Covering minima
and lattice point free convex bodies, {\it Annals of Math.}
{\bf 128} (1988), 577--602.

\item{[Ka]} {\smc J--M. Kantor}, Triangulations of integral
polytopes and Ehrhart polynomials, {\it Beitr\"age zur
Algebra und Geometrie}, {\bf 39} (1998), 205--218.

\item{[Lo]} {\smc L. Lov\'asz}, {\it An algorithmic theory of
numbers, graphs and convexity}, Regional Conference Series
in Applied Mathematics {\bf 50}, 1986.

\item{[McM]} {\smc P. McMullen}, Valuations and dissections,
in: {\it Handbook of convex geometry, ed. P. M. Gruber, J.
Wills}, North Holland, Amsterdam, 1988.

\enddocument

\bye